\documentclass{article}
\usepackage{graphicx, url, hyperref, geometry, color, xypic}
\input epsf.sty
\epsfverbosetrue 

\usepackage{latexsym}
\usepackage{amsfonts}
\usepackage{amssymb}

\newcommand{\R}{\mathbb{R}}
\newcommand{\Q}{\mathbb{Q}}
\newcommand{\Z}{\mathbb{Z}}
\newcommand{\C}{\mathbb{C}}
\newcommand{\IP}{\mathbb{P}}

\newcommand{\rs}{\mbox{$\widehat{\C}$}}

\def\MMM{{\mathcal M}}

\def\MMM{{\mathcal M}}
\def\WWW{{\mathcal W}}

\newtheorem{thm}{Theorem}[section]

\newtheorem{conj}[thm]{Conjecture}
\newtheorem{question}[thm]{Question}

\newcommand{\Imag}{\mbox{\rm Im}}                    

\newcommand{\Mod}{\mbox{\rm Mod}}    
\newcommand{\Teich}{\mbox{\rm Teich}}
\newcommand{\mod}{\mbox{\rm mod}}    
\newcommand{\union}{\cup}                        

\newcommand{\mtwo}[4]                            
{\mbox{$\left(\begin{array}{cc}                  
#1 & #2 \\
#3 & #4 
\end{array}
\right)$}}

\newcommand{\dettwo}[4]                          
{\mbox{$\left|\begin{array}{cc}                  
#1 & #2 \\
#3 & #4 
\end{array}
\right|$}}

\newcommand{\be}{\begin{enumerate}}
\newcommand{\eb}{\end{enumerate}}
\newcommand{\bi}{\begin{itemize}}
\newcommand{\ib}{\end{itemize}}
\newcommand{\bl}{\begin{list}}
\newcommand{\lb}{\end{list}}

\newcommand{\PSL}{\mbox{PSL}}

\newcommand{\cG}{\mathcal{G}}
\newcommand{\cT}{\mathcal{T}}
\newcommand{\cC}{\mathcal{C}}
\newcommand{\cA}{\mathcal{A}}
\newcommand{\cW}{\mathcal{W}}

\newcommand{\IR}{\mathbb{R}}

\newcommand{\pullback}{\tiny \stackrel{f}{\leftarrow}}
\newcommand{\wtgamma}{\widetilde{\gamma}}
\newcommand{\wtalpha}{\widetilde{\alpha}}

\newcommand{\BrMod}{\mathrm{BrMod}}
\newcommand{\Tw}{\mathrm{Tw}}
\begin{document}

\title{\bf On the pullback relation on curves induced by a Thurston map}
\author{Kevin M. Pilgrim\\Indiana University, Bloomington\\ \url{pilgrim@indiana.edu} }

\maketitle 

\abstract{Via taking connected components of preimages, a Thurston map $f: (S^2, P_f) \to (S^2, P_f)$ induces a pullback relation on the set of isotopy classes of curves in the complement of its postcritical set $P_f$.  We survey known results about the dynamics of this relation, and pose some questions.}

\section{Introduction}

An orientation-preserving branched covering $f: S^2 \to S^2$ of degree at least two is a \emph{Thurston map} if its \emph{postcritical set} $P_f = \union_{n>0}f^n(C_f)$ is finite, where $C_f$ is the finite set of branch (critical) points at which $f$ fails to be locally injective. 

A fundamental theorem in complex dynamics--Thurston's Characterization and Rigidity Theorem \cite{DH1}--asserts that apart from a well-known and ubiquitous set of counterexamples, the dynamics of rational Thurston maps is determined, up to holomorphic conjugacy, by its conjugacy-up-to-isotopy-relative-to-$P_f$ class.

Suppose $P \subset S^2$ is finite. The set of isotopy classes relative to $P$ of Thurston maps $f$ for which $P_f=P$ admits the structure of a countable semigroup under composition; we denote this by $\BrMod(S^2, P)$.  Pre- and post-composition with homeomorphisms fixing $P$ gives this semigroup the additional structure of a biset over the mapping class group $\Mod(S^2, P)$.   In this way, $\BrMod(S^2, P)$ may be fruitfully thought of as a generalization of the mapping class group.  This perspective is   useful in developing intuition for the range of potential behavior of and structure theory for Thurston maps.  

The mapping class group of a surface acts naturally on the countably infinite set of isotopy classes of curves on the surface. Even better, it acts on the associated curve complex; see \cite{MR624817}. It is natural to try to do something similar for Thurston maps.  Since the set $P_f$ contains the branch values of $f$, the restriction $f: S^2-f^{-1}(P_f) \to S^2-P_f$ is a covering map. It follows that a component $\wtgamma$ of the inverse image $f^{-1}(\gamma)$ of a simple closed curve $\gamma$ in $S^2-P_f$ is a simple closed curve in $S^2-f^{-1}(P_f)$. Since $P_f$ is forward-invariant, we have an inclusion $S^2-f^{-1}(P_f) \hookrightarrow S^2-P_f$,  so the curve $\wtgamma$ is again a simple closed curve in $S^2-P_f$. Abusing terminology, we'll call $\wtgamma$ a \emph{preimage} of $\gamma$, or sometimes say $\gamma$ \emph{lifts}, or \emph{pulls back}, to $\wtgamma$.  By lifting isotopies, we obtain a pullback relation $\pullback$ on the set of such simple closed curves $\cC$ up to isotopy.  The curve $\gamma$ might have several preimages, so we obtain an induced relation instead of a function. A preimage of an inessential curve is again inessential. Similarly, a preimage of a peripheral curve--one which is isotopic into any small neighborhood of a single point in $P_f$--is either again peripheral, or is inessential. We call inessential and peripheral curves \emph{trivial}, and note that the set of trivial curves is invariant under the pullback relation. 

When $\#P_f=4$, the pullback relation induces--almost--a function on the set of nontrivial curves. On the one hand, distinct nontrivial curves in this case must intersect. On the other hand, distinct components of $f^{-1}(\gamma)$ are in general disjoint.  It follows that there can be at most one class of nontrivial preimage, and we almost get a function in this case. Why ``almost''? Typical examples have the property that for some curve, each of its preimages are trivial. So while the mapping class group acts naturally on e.g. the infinite diameter curve complex, it is less clear how to construct a nice complex related to curves on which a Thurston map acts via pullback. This relative lack of preserved structure makes answering even basic questions challenging. 

This survey presents some known results about the dynamical behavior of taking iterated preimages of curves under a given Thurston map. It assumes the basic vocabulary related to Thurston maps from \cite{DH1}, and the reader may find \cite{kmp:kps} also useful for more detailed explanations and references to some terminology encountered along the way. 

Here are some highlights, to convince you that this is interesting.  When $f(z)=z^2+c$ is the so-called \emph{Douady Rabbit} polynomial, where $c$ is chosen so that $\Im(c)>0$ and the origin has period $3$, under iteration every curve pulls back to either a trivial curve, or to the prominent 3-cycle.  But when $f(z)=z^2+i$, every curve pulls back eventually to a trivial curve. See \cite{kmp:tw:published}. In these two cases, we see that there is a \emph{finite global attractor} for the pullback relation. Among obstructed Thurston maps, though, it is easy to manufacture examples with wandering curves and infinitely many fixed curves. A basic conjecture is 
\begin{conj}
\label{conj:fga}
If $f$ is rational and not a flexible Latt\`es example, then the pullback relation on curves has a finite global attractor.
\end{conj}
Here, informally, is the source of the tension. Represent a curve class  by the unique geodesic $\gamma$ in the complement of $P_f$ equipped with its hyperbolic metric. Pulling back and lifting the metric to $\rs - f^{-1}P_f$, the lifted curve $\wtgamma$ may unwind and become up to $\deg(f)$ times as long as $\gamma$. But when including the curve $\wtgamma$ back into $\rs-P_f$, the Schwarz lemma implies the length of $\wtgamma$ shrinks. It is unclear which force--lengthening or shortening--dominates in the long run. And since there exist expanding non-rational examples with wandering curves, the exact mechanism that would imply the conjecture remains mysterious. 

\subsection*{Acknowledgements} This work was partially supported by Simons Foundation collaboration grants 4429419 and 615022.

\section{Conventions and notation}
Throughout, $f$ denotes a Thurston map, $P$ its postcritical set, and $d$ its degree. To avoid repeated mention of special cases, unless otherwise stated, $f$ has hyperbolic orbifold and $\#P \geq 4$.  We denote by 
\bi
\item $\simeq$ the equivalence relation of isotopy-through-Thurston-maps-with-fixed-postcritical-set $P$;
\item $\cC$, the countably infinite set of istopy classes of unoriented, essential, simple, nonperipheral curves in $S^2-P$ (we will often call such elements simply ``curves'', abusing terminology); on $\cC$ we have the \emph{geometric intersection number} $\iota(\alpha, \beta)$ which counts the minimum number of intersection points among representatives; 
\item $o$, the union of the $\#P+1$ isotopy classes of unoriented, simple, closed, peripheral and inessential curves in $S^2-P$, i.e. the trivial ones;
\item $\overline{\cC}:=\cC \union \{o\}$;
\item $\pullback$, the pullback relation on $\overline{\cC}$ induced by $\gamma \mapsto \delta \subset f^{-1}(\gamma)$, where $[\gamma]\in\overline{\cC}$ and $\delta$ is a component of $f^{-1}(\gamma)$;
\item $\overline{\cA}$ and $\cA$, the set of curves contained  in cycles of $\pullback$ in $\overline{\cC}$  and $\cC$, respectively;
\item $\cW \subset \cC$, the set of ``wandering'' curves $\gamma_0$, namely, those for which there is an infinite sequence $\gamma_n, n \geq 0$, of distinct nontrivial curves satisfying $\gamma_n \pullback \gamma_{n+1}, n\geq 0$;
\item the relation $\pullback$ has a \emph{finite global attractor} if $\cW$ is empty and $\cA$ is finite; 
\item $\Teich(S^2, P)$, the Teichm\"uller space of the sphere marked at the set $P$;
\item $\sigma_f: \Teich(S^2,P) \to \Teich(S^2,P)$, the holomorphic self-map obtained by pulling back complex structures; it is the lift to the universal cover of an algebraic \emph{correspondence on moduli space} $X\circ Y^{-1}$, where $Y$ is a finite cover and $X$ is holomorphic. See Figure 1. 
\begin{figure}[h]
\label{figure:w} 
\[
\xymatrix{  & \Teich(S^2, P) \ar[dd]_{\pi}\ar[rr]^{\sigma_f} \ar[dr]^{\omega_f} & &
\Teich(S^2, P) \ar[dd]^{\pi} \\
&&\WWW_f\ar[dl]_{Y_f} \ar[dr]^{X_f} &\\
& \MMM_P & & \MMM_P}
\]
\caption{\sf The fundamental diagram.}\label{fund}
\end{figure}
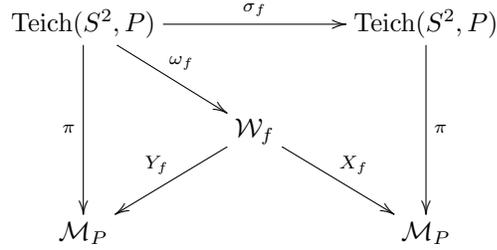
\ib 

\section{Non-dynamical properties of $\pullback$.} 

A Thurston map $f: (S^2, P) \to (S^2, P)$ may factor as a composition of maps of pairs  
\[(S^2, P) \stackrel{f_1}{\to} (S^2, C) \stackrel{f_2}{\to} (S^2, P)\] where each $f_i $ is admissible in the sense that its set of branch values is contained in the distinguished subset appearing in its codomain. This motivates studying properties of so-called admissible branched covers $f: (S^2, A) \to (S^2, B)$ where domain and codomain are no longer identified; this perspective was introduced by S. Koch.  Instead of a pullback self-relation on curves $\cC$, we have a pullback relation $\cC_B \pullback \cC_A$ from classes of curves in $S^2-B$ to classes in $S^2-A$. 

\subsection{Known general results}
Thinking non-dynamically first, we have the following known results about the pullback relation $\pullback$. 
\be
\item When $\#A=\#B$, each nonempty fiber is dense in the Thurston boundary; in particular, each nontrivial fiber is infinite \cite{kmp:kps}.  Here, by the fiber over $\beta$, we mean $\{ \alpha : \alpha \pullback \beta\}$.  
\item The relation $\pullback$ can be trivial in the sense that the only pairs are of the form  $\gamma \pullback \delta$ where $\delta$ is trivial.  Equivalently, $\sigma_f$ is constant.  See \cite{kmp:kps}, correcting an argument appearing originally in \cite{bekp}. 
\item The relation $\pullback$ satisfies a Lipschitz-type inequality related to intersection numbers: $\iota(\widetilde{\alpha}, \widetilde{\beta}) \leq d\cdot \iota(\alpha, \beta)$ whenever $\alpha \pullback \wtalpha, \beta \pullback \widetilde{\beta}$.

The study of the interaction between intersection numbers and the geometry of $\sigma_f$ seems to be just beginning.  Implicit use of such interactions appears in the analysis of the family of so-called Nearly Euclidean Thurston (NET) maps by W. Floyd, W. Parry, and this author; see \cite{cfpp:net}, \cite{Floyd:2018aa}. Parry develops this intersection theory further in \cite{parry:netmap:slopefunctions}.   A careful study is also applied in \cite{bonk:hlushchanka:iseli:curves} to give new methods for analyzing the effect of certain surgery operations on possible obstructions. 

\item Proper multicurves are in natural bijective correspondence with boundary strata in the augmented Teichm\"uller space, which by a theorem of Masur is known to be the completion of Teichm\"uller space in the Weil-Petersson metric \cite{masur:duke:1976}. A result of Selinger \cite{selinger:pullback} shows that $\sigma_f: \Teich(S^2,A) \to \Teich(S^2,B)$ extends to the this completion, sending the stratum corresponding to a multicurve $\Gamma$ to the stratum corresponding to the multicurve $f^{-1}(\Gamma)$. It follows that analytical tools for studying $\sigma_f$ can be used to study properties of the combinatorial relation $\pullback$ \cite{kmp:tw:published}, \cite{kmp:kps}.  There is thus a rich interplay between the analytic and algebro-geometric properties of the correspondence on moduli space, and the combinatorial properties of the pullback relation; see recent work of R. Ramadas, for example.  

\item Associated to a proper multicurve $\Gamma$ is the free abelian group $\Tw(\Gamma)$ of products of powers of Dehn twists about the curves in $\Gamma$. The pullback function can be encoded using the associated induced virtual endomorphism on the mapping class group $\phi_f: \Mod(S^2, P_f) \dashrightarrow \Mod(S^2, P_f)$ induced by lifting. It follows that group-theoretic tools can also be used to study properties of the pullback relation on curves; see \cite{kmp:tw:published} and \cite{lodge:kelsey:quadratic}.  
\eb

\begin{question}
\label{q:surjective}
If the pullback relation $\pullback$ is not trivial, must it be surjective?
\end{question}

It seems very likely that the answer is no, for the following reason. The Composition Trick, introduced in the next subsection, should allow one to build examples where the image of $\sigma_f$ has positive dimension and codimension, so that its image misses many strata. 

\subsection{Mechanisms for triviality of $\pullback$}
There seem to be three or four mechanisms via which $\pullback$ can be trivial. 

\be
\item {\bf Composition trick.}  The map $f:(S^2,A) \to (S^2,B)$ may factor through $(S^2, C)$ with $\#C=3$ (C. McMullen, \cite{bekp}). Such maps have the property that $\sigma_f$ and $\pullback$ are trivial. 

\item {\bf NET maps.}   A. Saenz \cite{saenz:thesis} found an example of a Thurston map $f$ for which $\sigma_f$ is constant but for which $f$ does not decompose as in the Composition Trick. Here is his example, from a different point of view. 

Let $E$ be an elliptic curve over $\C$. Regarded as an abelian group, there are 8 distinct points of order $3$; under the involution $z \mapsto -z$ these 8 points descend to a set $A$ of 4 points on $\IP^1=E/\pm 1$ whose cross-ratio is, miraculously, constant as $E$ varies. Now take $E$ to be the square torus and let $f: \IP^1 \to \IP^1$ be the degree 9 flexible Latt\`es map induced by the tripling map on $E$, and let $B=f(A)$. As $E$ varies, the conformal shape of $B$ varies, but that of $A$ does not. Thus $\sigma_f$ is constant and so $\pullback$ is trivial. One can see this triviality directly by observing that the action of $\PSL_2(\Z)$ is transitive on curves (since it acts transitively on extended rationals regarded as slopes), that $A$ is invariant under this action (since points of order 3 are invariant under the induced group-theoretic automorphisms), and that the horizontal curve has all preimages inessential or peripheral (as drawing a single easy picture shows).

\item {\bf Sporadic examples.} Let $f$ be the unique (up to pre- and post-composition by independent automorphisms) degree four rational map with three double critical points $(c_1, c_2, c_3)$ mapping to necessarily distinct critical values $(v_1, v_2, v_3)$. Choose $w$ a point distinct from the $v_i$'s, let $B=\{v_1, v_2, v_3, w\}$ and $A=R^{-1}(w)=\{z_1, z_2, z_3, z_4\}$. Then the $j$-invariant (obtained from the cross-ratio by applying a certain degree six rational function) of the $z_i$'s is constant in $w$, whence $\sigma_f$ is constant and so $\pullback$ is trivial. To see this, note that as $w$ approaches some $v_i$, three of the $z_i$'s converge to $c_i$, and the remaining one converge to some other point, call it $c_i'$. Normalizing so $c_i=0$ and $c_i'=\infty$ and scaling via multiplication with a nonzero complex constant shows that the conformal shape of the fiber $R^{-1}(w)$ converges to that of the cube roots of unity together with the point at infinity.  Thus the $j$-invariant of $R^{-1}(w)$ is a bounded holomorphic function, hence constant. 

\item {\bf Combinations of the above.}  
\eb
\begin{question}
\label{q:primedeg}
Do there exist examples $f$ with $\sigma_f$ constant and $\deg(f)$ prime?
\end{question} 

Cui G. has thought about the general case, see \cite{cui-constant-sigma-slides}. It is natural to look for the simplest such examples. By cutting along a maximal multicurve, one may restrict to the case $\#A=\#B=4$; let's call these ``minimal''. It is natural to look for examples which do not factor as in the Composition Trick; let's call these ``primitive''. 

\begin{question}
\label{q:primitive_trivial}
What are the minimal primitive branched covers $f: (S^2,A) \to (S^2,B)$ for which $\pullback$ is trivial?
\end{question}

\subsection{Computation of $\pullback$}

Though the set of curves $\cC$ is complicated, it is conveniently described by a variety of coordinate systems. For example, Dehn-Thurston coordinates record intersection numbers of curves with edges in a fixed triangulation with vertex set $P$ \cite{MR3053012}.  Train tracks and measured foliations give other methods. But expressing the pullback relation in these coordinates can be very complicated. The effect on Dehn-Thurston coordinates of pulling back from $S^2-P_f$ to $S^2-f^{-1}(P_f)$ is indeed easy to compute, since one can just lift intersection numbers. But the expression for the induced ``erasing map'' is hard to write down in closed form and leads to continued-fraction-like cases. 

When $\#P_f=4$, though, the set of curves $\cC$ can be encoded by ``slopes'' in the extended rationals $\Q \union \{1/0\}$, the pullback relation $\pullback$ is a almost a function, and things are a bit easier but are still quite complicated. 

A Thurston map is an NET map if $\#P_f=4$ and each critical point has local degree two; see \cite{cfpp:net}.  If $f$ is an NET map, there is an algorithm that computes the image of a slope under pullback. This can be done easily by hand, and has been implemented. NET maps can be easily encoded by combinatorial input. W. Parry has written a computer program that implements this algorithm. The website \url{http://intranet.math.vt.edu/netmaps/}, maintained by W. Floyd, contains a database of tens of thousands of examples.  For NET maps, it appears that this ability to calculate the pullback map on curves (and related invariants, such as the degree by which preimages map, and how many preimages there are) leads to an effective algorithm for determining whether a given example is, or is not, equivalent to a rational map. 

When $\#P_f=4$ and $f$ is the subdivision map of a subdivision rule of the square pillowcase (like in Figure \ref{fig:blown_up_lattes}), W. Parry has written a program for computing the image of a slope under pullback (personal communication). 

\begin{question}
\label{q:computer_in_higher}
Are there any interesting settings in which one can effectively compute $\pullback$ when $\#P_f \geq 5$? 
\end{question}

\subsection{When each curve has a nontrivial preimage}

The example studied by Lodge \cite{lodge:thesis:ecgd} is, nondynamically speaking, the unique generic cubic, in the following sense. We have $f: (S^2, A) \to (S^2, B)$ where $A$ is the set of four simple critical points and $B=f(A)$. Each nontrivial curve has exactly one nontrivial curve in its preimage.

\begin{question}
\label{q:generic_is_proper}
Suppose $f: (S^2,A) \to (S^2,B)$ has the property that $A$ consists of $2d-2$ simple critical points, $B=f(A)$, and $f|B$ is injective, so that $B$ consists of the $2d-2$ critical values. Does each nontrivial curve in $S^2-B$ have a nontrivial preimage?
\end{question}

Up to pre- and post-composition with homeomorphisms there is a unique such map \cite{berstein:edmonds:construction}.  If each $\gamma$ has a nontrivial preimage, so does $h(\gamma)$, where $h: (S^2,P_f) \to (S^2, P_f)$ is a homeomorphism that lifts under $f$. Since the set of such $h$ is a finite-index subgroup of $\Mod(S^2, P_f)$, checking this is a finite computation. 

In the case of four postcritical points, we have the following result from \cite{Floyd:2018aa}.
\begin{thm}
\label{thm:trivial_or_surjective}If $\#P=4$, the pullback function on curves is either surjective, or trivial.
\end{thm} 
The key insight: when $\#P=4$, looking at the correspondence on moduli space, the space $\cW_f$ a Riemann surface whose set of ends consists of finitely many cusps, the map $X: \cW_f \to \MMM_P$ is holomorphic, and $\MMM_P$ is the triply-punctured sphere; thus if $X$ is nonconstant, then each cusp of $\MMM_P$ is the image of cusp of $\cW_f$.

\section{Dynamical properties} 

\subsection{General properties} We first discuss in general some examples and known results about the possible dynamical behavior of $\pullback$.

\be
\item {\bf Example:} \emph{Every curve iterates to the trivial curve.} This happens for $z^2+i$.  Here is one way to see this.  Examining the possibilities for how the bounded region enclosed by a curve meets the finite postcritical set $\{i, i-1, -i\}$, one sees that a curve must eventually become trivial unless it surrounds both $-i$ and $i-1$. For this type of curve $\alpha$, there is at most one nontrivial curve $\beta$ with $\alpha \pullback \beta$ and $\beta$ a curve of the same type.  Moreover, $\deg(\alpha \pullback \beta) = 1$.  Equipping the complement of the postcritical set with the hyperbolic metric, the Schwarz Lemma shows that the length of a geodesic representative of $\beta$ is strictly shorter than that of $\alpha$.   Iterating this process, it follows that such a curve cannot be periodic under $\pullback$: points in its orbit cannot get too complicated, since otherwise they would have to get long, so they must eventually become a different type of curve and thus become trivial upon further iteration. 

The ``airplane'' quadratic polynomial $f(z)=z^2+c$, with the origin periodic of period 3 and $\Imag(c)=0$, is another such example \cite{lodge:kelsey:quadratic}.  

\item 
\begin{question}
\label{q:nilpotent}
Does there exist an example of a Thurston map $f$ for which the pullback relation induced by $f$ is nontrivial but that induced by some iterate $f^n$ is trivial? 
\end{question}

\item {\bf Theorem.}  If $f$ is rational and non-Latt\`es, then there are only finitely many fixed proper multicurves for which $f^{-1}\Gamma = \Gamma$; see \cite[Thm. 1.5]{kmp:tw:published}. The proof uses the decomposition theory.  

\item {\bf Conjecture 1.1:} \emph{If $f$ is rational and not a flexible Latt\`es example then the pullback relation on curves has a finite global attractor.}  

There is partial progress on this conjecture for special families of Thurston maps.
\be
\item Kelsey and Lodge \cite{lodge:kelsey:quadratic} verify this for all quadratic non-Latt\`es maps with four postcritical points.

\item Hlushchanka \cite{Hlushchanka:2019aa} verifies this for critically fixed rational maps. As he shows, each such map is obtained by blowing up, in the sense of \cite{kmp-tan:blow}, edges of a planar multigraph $\cG$. If $f$ is a critically fixed rational map, $\alpha$ an edge in the planar connected multigraph $\cG$ that describes $f$ via the blowing up construction \cite{kmp:fixed}, $\gamma \pullback \wtgamma$, and $\iota(\gamma, \alpha) > 0$, then it is easy to see that unless $\gamma$ is homotopic to the boundary of a regular neighborhood of an edge-path homeomorphic to an embedded arc in $\cG$ and $\iota(\gamma,\alpha)=1$, we must have $\iota(\wtgamma, \alpha)<\iota(\gamma,\alpha)$.    The global attractor $\cA$ consists of those $\gamma$ for which $\iota(\gamma,\alpha) \leq 1$ for all edges $\alpha \in \cG$. 

\item If the virtual endomorphism $\phi_f$ on the mapping class group is contracting, then $\pullback$ has a finite global attractor \cite{kmp:kps}. Nekrashevych \cite[Thm. 7.2]{vn:models:MR3199801} establishes this contraction property in the case for hyperbolic polynomials.

\item Belk, Lannier, Margalit, and Winarski \cite{Belk:2019aa} associate to a topological polynomial $f$ (a Thurston map with a fixed point at infinity which is fully ramified) a simplicial complex $\cT$ whose vertices are planar trees, and a simplicial map $\lambda_f: \cT \to \cT$ induced by lifting. For general complex polynomials, the associated Hubbard trees are fixed vertices, and the uniquness of the Hubbard tree for iterates of $f$ leads to a contraction property of $\lambda_f$ that implies the existence of a finite global attractor for the pullback relation on curves. 

\begin{question}
\label{question:simplicial_action}
For a general Thurston map, is there an associated natural action on a contractible simplicial complex?
\end{question}

\item If the correspondence on moduli space (in the direction of $\sigma_f$) has a nonempty invariant compact subset, then $\phi_f$ is contracting, so there is a finite global attractor.  If moduli space admits an incomplete metric which is (i) uniformly contracted by $\sigma_f$, and (ii) whose completion is homeomorphic to that of the WP metric, then the trivial curve is a finite global attractor \cite[Thm. 7.2]{kmp:kps}. The latter occurs for $f(z)=z^2+i$; the correspondence on moduli space is the inverse of a Latt\`es map with three postcritical points and Julia set the whole sphere, which expands the Euclidean orbifold metric. 

\item Intersection theory provides some insight. Consider the linear map $\lambda_f: \IR[\cC] \to \IR[\cC]$ defined on basis vectors by $\gamma \mapsto \sum_{\gamma \pullback \delta}\frac{1}{d_\delta}\cdot \delta$ where $\deg(\gamma \pullback \delta)=d_\delta$. Then one can show that the Dehn-Thurston coordinates of any orbit of $\lambda_f$ starting at any curve $\gamma$ tends to zero. For otherwise, one has either an obstruction, ruled out by rationality, or a wandering curve with a uniform lower bound in the corresponding coordinates and hence on the corresponding moduli of path families.  But if $\alpha, \beta$ are any two curves in $\rs - P$ then $\mod(\alpha)\mod(\beta) \lesssim 1/\iota(\alpha, \beta)^2$ and so this is impossible. 
\eb
\eb

\subsection{Bounds on the size of the attractor}  Since up to conjugacy there are only finitely many non-flexible Latt\`es rational maps with a given degree $d$ and size $\#P$ of postcritical set, the cardinality of the finite global attractor $\cA$, if one exists, must be bounded by some constant depending on $d$ and $\#P$.   I know very little about the behavior of this function.  

\be
\item Certainly $\# \cA$ can be large if $\#P$ is large: for $n \geq 2$ the ``$1/n$-rabbit'' quadratic polynomial will have an $n$-cycle of curves. Other examples can be constructed by perturbing flexible Latt\`es examples. One can find hyperbolic sets consisting of invariant curves which are stable under perturbation; a result of X. Buff and T. Gauthier \cite[Cor. 3]{MR3203689} implies that such maps are limits of sequences of hyperbolic Thurston maps with the maximum number $2d-2$ attracting cycles. This observation gives robust constructions for examples with many fixed curves.

\item In composite degrees, $\#\cA$ can be small (say zero), by taking e.g. examples with $\sigma_f$ constant. Using McMullen's compositional trick and Belyi functions one can easily build both hyperbolic rational maps and rational maps with Julia set the whole sphere having the property that $\sigma_f$ constant and $\#P_f$ arbitrarily large. 

\item Results of G. Kelsey and R. Lodge \cite{lodge:kelsey:quadratic} show that for quadratic rational maps $f$ with $\#P_f=4$, we have $\#\cA \leq 4$.   

The bound might be explained as follows.  The map $f$ corresponds (not quite bijectively) to a repelling fixed-point $p$ of a correspondence $g=Y\circ X^{-1}$ on moduli space. In the nonexceptional cases, this is actually a rational map $g: \IP^1 \to \IP^1$. There appears to be a natural bijection between invariant (multi)curves for $f$ and periodic internal rays joining points in periodic superattracting cycles of $g$ (these lie at infinity in moduli space) to $p$.  I've confirmed this also for critically fixed polynomials with three finite critical points.  

\begin{question}
\label{q:inv_mc_and_internal_rays}
If $f$ is rational and $\Gamma$ is an invariant or periodic multicurve, does there exist a stable manifold connecting the unique fixed-point of $\sigma_f$ in $\Teich(\rs, P)$ to a fixed-point or periodic point in a boundary stratum corresponding to $\Gamma$? 
\end{question}

For quadratics with four postcritical points, the analysis of \cite{lodge:kelsey:quadratic} seems to confirm the intuition that periodic curves in the dynamical plane are related to accesses landing at the associated fixed-point that are periodic under the correspondence, in this case the inverse of a critically finite rational map with three postcritical points.  But in higher degrees with $\#P_f=4$, the correspondence need not the inverse of such a map, and the situation is more complicated; see e.g. \cite{lodge:thesis:ecgd}.  

In higher dimensions, another first natural example to try is the case of $f$ a critically fixed polynomial with four finite simple critical points. One would need to show the existence of internal rays in two complex dimensions.  This example is beautifully symmetric and possesses many invariant lines that might make the problem more tractable. See \cite[section 3]{bek:bottcher}.

\eb

\subsection{Examples with symmetries}  Maps with nontrivial symmetries provide a source of non-rational examples without a finite global attractor. We denote by $\Mod(f)=\{h : hf\simeq fh\; \mbox{rel}\; P_f\}$; here $\simeq$ denotes isotopy. We recall four facts:
\be
\item The pure mapping class group $P\Mod(S^2, P)$ has no elements of finite order, so neither does $P\Mod(f)$. 

\item If $f$ is rational, $P\Mod(f)$ is trivial, unless $f$ is a flexible Latt\`es example, in which case it is isomorphic to the free group on two generators.

\item Suppose $f$ has an obstruction $\Gamma=(\gamma_1, \ldots, \gamma_m)$ with the property that the corresponding linear map $f_\Gamma: \R[\Gamma] \to \R[\Gamma]$ has $1$ among its eigenvalues, with a corresponding nonnegative integer eigenvector $(a_1, \ldots, a_m)$. Let $T_i$ denote the Dehn twist about $\gamma_i$. Then some power of $T_1^{a_1}\cdots T_m^{a_m}$ gives an element of $\Mod(f)$ \cite{kmp:cds}.

\item Thurston maps are like mapping classes.  If $f$ is obstructed, there is a canonical decomposition by cutting along a certain invariant multicurve \cite{kmp:canonical}, \cite{selinger:pullback}. The ``pieces'' might contain cycles of degree one: mapping class elements, each with its own centralizer.  The fact that the decomposition is canonical means that the centralizers of the pieces will embed into $\Mod(f)$.  Using this idea one can create examples of Thurston maps with a variety of prescribed behaviors. For example, if $f$ is the identity on some sufficiently large piece, then clearly $\cA$ contains infinitely many fixed curves. If $f$ is a pseudo-Anosov map on some other sufficiently large piece, then there are infinitely many distinct orbits of wandering curves.  
\item L. Bartholdi and D. Dudko give an explicit example of $f$ with $\Mod(f)$ infinitely generated \cite{bartholdi:dudko:bc:ii}. 

\eb
\subsection{Maps with the same fundamental invariants} 

As motivation, recall that if $L$ is a flexible Latt\`es example with postcritical set $P$, then $\sigma_L$ acts as the identity, and the pullback relation on curves is the identity function. So if $f$ is now an arbitrary Thurston map with with the same postcritical set $P$, then $f\circ L \simeq L \circ f$, and $\sigma_{L\circ f}=\sigma_f$ and the pullback relation on curves for $f$ and $L\circ f$ are the same.


\begin{question}
\label{q:sameinv}
When do two Thurston maps have the same pullback relation on curves? 
\end{question}

\subsection{Expanding vs. nonexpanding maps}  The examples in \S 4.3(4) are  not isotopic to expanding maps. There exist Levy cycles--cycles in which each curve maps by degree 1--and these are obstructions to the existence of an expanding representative. However, there exist expanding maps without finite global attractors. 

Blow up the 2x2 Latt\`es doubling example along the middle upper vertical edge to get a Thurston map $f$; see Figure 2.
\begin{figure}[h]
\begin{center}
\includegraphics[width=3in]{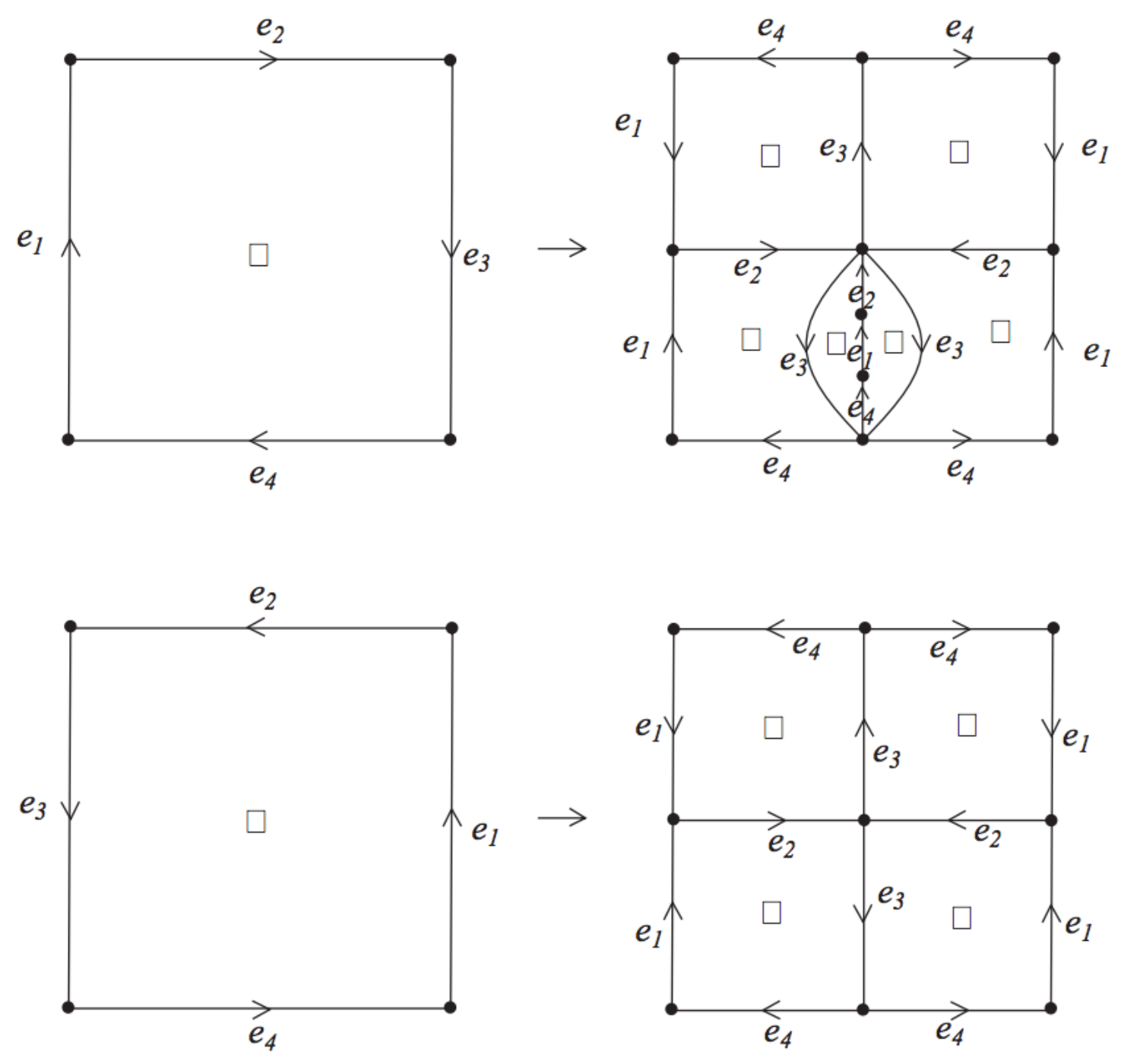}
\caption{The codomain is the union of the two squares at left along their boundaries as indicated to form a square ``pillowcase''. Each square at right is identified with the square just to its left by a translation, so that the pillowcases are identified. The figure shows a cell structure in domain and codomain.  The map $f$ goes in the opposite direction to the indicated arrow and defines a cellular degree 5 map from the pillowcase to itself. The four corners of the pillowcase form the postcritical set. \emph{Figure by W. Floyd.}  }
\end{center}
\end{figure}
The ``rim'' of the square pillowcase--the common boundary of the two squares at left-- is an invariant Jordan curve containing $P_f$: that is, $f$ is the subdivision map of a finite subdivision rule.   This map satisfies the combinatorial expansion properties of both Bonk-Meyer \cite{MR3727134} and Cannon-Floyd-Parry \cite{cfp:fsr}.  Appealing to either one of these works, we conclude there is a map isotopic to this example in which the diameters of the tiles goes to zero upon iterated subdivision. This implies that this example has no Levy cycles. Appealing to \cite{MR3852445}, we conclude there is such an example that is expanding with respect to a complete length metric. 

The vertical curve is an obstruction with multiplier 1, and the horizontal curve is invariant. Let $T$ be the Dehn twist about this vertical curve, so that $fT=Tf$. This immediately implies (i) $\cA_f$ is infinite, since the orbit of the horizontal curve under $T$ will consist of $f$-invariant curves, and (ii) if we put $g=Tf$,  then the horizontal curve wanders. To see this, note that $g^{-n}(\gamma)=(Tf)^{-n}(\gamma)=T^{-n}f^{-n}(\gamma)=T^{-n}\gamma$ as required. I do not know if $Tf$ is isotopic to an expanding map. Letting $L$ denote the flexible Latt\`es map induced by doubling on the torus, however, I would expect that $TfL^N$ is expanding for sufficiently large $N$. This should show the existence of expanding maps with wandering curves. 

\newpage

\begin{question}
\label{q:assymetric_with_wandering}
Suppose $\Mod(f)$ is trivial.  Could there exist infinitely many periodic curves? Could there exist wandering curves? 
\end{question}

\bibliographystyle{/Users/pilgrim/Dropbox/FMEO/math}
\bibliography{/Users/pilgrim/Dropbox/FMEO/kmprefs.bib}

\newcommand{\etalchar}[1]{$^{#1}$}
\def\cprime{$'$}
\begin{thebibliography}{CGN{\etalchar{+}}}

\bibitem[BD]{bartholdi:dudko:bc:ii}
Laurent Bartholdi and Dzmitry Dudko.
\newblock {Algorithmic aspects of branched coverings II/V: sphere bisets and
  their decompositions.}
\newblock \url{https://arxiv.org/pdf/1603.04059.pdf}.

\bibitem[BE]{berstein:edmonds:construction}
Israel Berstein and Allan~L. Edmonds.
\newblock {On the construction of branched coverings of low-dimensional
  manifolds}.
\newblock {\em Trans. Amer. Math. Soc.} {\bf 247}(1979), 87--124.

\bibitem[BM]{MR3727134}
Mario Bonk and Daniel Meyer.
\newblock {\em Expanding {T}hurston maps}, volume 225 of {\em Mathematical
  Surveys and Monographs}.
\newblock American Mathematical Society, Providence, RI, 2017.

\bibitem[BEKP]{bekp}
Xavier Buff, Adam Epstein, Sarah Koch, and Kevin Pilgrim.
\newblock {On {T}hurston's pullback map}.
\newblock In {\em Complex dynamics}, pages 561--583. A K Peters, Wellesley, MA,
  2009.

\bibitem[BG]{MR3203689}
Xavier Buff and Thomas Gauthier.
\newblock {Perturbations of flexible {L}att\`es maps}.
\newblock {\em Bull. Soc. Math. France} {\bf 141}(2013), 603--614.

\bibitem[CFP]{cfp:fsr}
J.~W. Cannon, W.~J. Floyd, and W.~R. Parry.
\newblock {Finite subdivision rules}.
\newblock {\em Conform. Geom. Dyn.} {\bf 5}(2001), 153--196 (electronic).

\bibitem[CFPP]{cfpp:net}
J.~W. Cannon, W.~J. Floyd, W.~R. Parry, and K.~M. Pilgrim.
\newblock {Nearly {E}uclidean {T}hurston maps}.
\newblock {\em Conform. Geom. Dyn.} {\bf 16}(2012), 209--255.

\bibitem[CGN{\etalchar{+}}]{kmp:fixed}
Kristin Cordwell, Selina Gilbertson, Nicholas Nuechterlein, Kevin~M. Pilgrim,
  and Samantha Pinella.
\newblock {On the classification of critically fixed rational maps}.
\newblock {\em Conform. Geom. Dyn.} {\bf 19}(2015), 51--94.

\bibitem[Cui]{cui-constant-sigma-slides}
Guizhen Cui.
\newblock {Rational maps with (constant) pullback map}.
\newblock
  \url{http://www.math.ac.cn/kyry/cgz/201501/W020150121410977178324.pdf}.

\bibitem[DH]{DH1}
A.~Douady and John Hubbard.
\newblock {A Proof of {T}hurston's Topological Characterization of Rational
  Functions}.
\newblock {\em Acta. Math.} {\bf 171}(1993), 263--297.

\bibitem[FKK{\etalchar{+}}]{kmp:origami}
William Floyd, Gregory Kelsey, Sarah Koch, Russell Lodge, Walter Parry,
  Kevin~M. Pilgrim, and Edgar Saenz.
\newblock {Origami, affine maps, and complex dynamics}.
\newblock {\em Arnold Mathematical Journal} {\bf 3}(2016), 365--395.

\bibitem[FPP]{Floyd:2017fj}
William Floyd, Walter Parry, and Kevin~M. Pilgrim.
\newblock {Presentations of NET maps}.
\newblock (01 2017).

\bibitem[KL]{lodge:kelsey:quadratic}
Gregory Kelsey and Russell Lodge.
\newblock {Quadratic Thurston maps with few postcritical points,
  \url{https://arxiv.org/abs/1704.03929}}.
\newblock (04 2017).

\bibitem[KPS]{kmp:kps}
Sarah Koch, Kevin~M. Pilgrim, and Nikita Selinger.
\newblock {Pullback invariants of {T}hurston maps}.
\newblock {\em Trans. Amer. Math. Soc.} {\bf 368}(2016), 4621--4655.

\bibitem[Lod]{lodge:thesis:ecgd}
Russell Lodge.
\newblock {Boundary values of the {T}hurston pullback map}.
\newblock {\em Conform. Geom. Dyn.} {\bf 17}(2013), 77--118.

\bibitem[Mal]{saenz:thesis}
Edgar Arturo~Saenz Maldonado.
\newblock {On Nearly Euclidean Thurston Maps}.
\newblock {\em PhD Thesis, Virginia Polytechnic University} (2012).

\bibitem[Pil1]{kmp:cds}
Kevin~M. Pilgrim.
\newblock {\em Combinations of complex dynamical systems}, volume 1827 of {\em
  Lecture Notes in Mathematics}.
\newblock Springer-Verlag, Berlin, 2003.

\bibitem[Pil2]{kmp:tw:published}
Kevin~M. Pilgrim.
\newblock {An algebraic formulation of {T}hurston's characterization of
  rational functions}.
\newblock {\em Ann. Fac. Sci. Toulouse Math. (6)} {\bf 21}(2012), 1033--1068.

\bibitem[Sel]{selinger:pullback}
Nikita Selinger.
\newblock {Thurston's pullback map on the augmented {T}eichm\"uller space and
  applications}.
\newblock {\em Invent. Math.} {\bf 189}(2012), 111--142.

\end{thebibliography}


\newcommand{\etalchar}[1]{$^{#1}$}
\def\cprime{$'$}
\begin{thebibliography}{BLMW}

\bibitem[BD1]{MR3852445}
Laurent Bartholdi and Dzmitry Dudko.
\newblock {Algorithmic aspects of branched coverings {IV}/{V}. {E}xpanding
  maps}.
\newblock {\em Trans. Amer. Math. Soc.} {\bf 370}(2018), 7679--7714.

\bibitem[BD2]{bartholdi:dudko:bc:ii}
Laurent Bartholdi and Dzmitry Dudko.
\newblock {Algorithmic aspects of branched coverings II/V: sphere bisets and
  their decompositions.}
\newblock \url{https://arxiv.org/pdf/1603.04059.pdf}.

\bibitem[BLMW]{Belk:2019aa}
James Belk, Justin Lanier, Dan Margalit, and Rebecca~R. Winarski.
\newblock {Recognizing topological polynomials by lifting trees,
  \url{https://arxiv.org/pdf/1906.07680.pdf}}, 06 2019.

\bibitem[BE]{berstein:edmonds:construction}
Israel Berstein and Allan~L. Edmonds.
\newblock {On the construction of branched coverings of low-dimensional
  manifolds}.
\newblock {\em Trans. Amer. Math. Soc.} {\bf 247}(1979), 87--124.

\bibitem[BM]{MR3727134}
Mario Bonk and Daniel Meyer.
\newblock {\em Expanding {T}hurston maps}, volume 225 of {\em Mathematical
  Surveys and Monographs}.
\newblock American Mathematical Society, Providence, RI, 2017.

\bibitem[BEKP]{bekp}
Xavier Buff, Adam Epstein, Sarah Koch, and Kevin Pilgrim.
\newblock {On {T}hurston's pullback map}.
\newblock In {\em Complex dynamics}, pages 561--583. A K Peters, Wellesley, MA,
  2009.

\bibitem[BEK]{bek:bottcher}
Xavier Buff, Adam~L. Epstein, and Sarah Koch.
\newblock {B{\"o}ttcher coordinates}.
\newblock {\em Indiana Univ. Math. J.} {\bf 61}(2012), 1765--1799.

\bibitem[BG]{MR3203689}
Xavier Buff and Thomas Gauthier.
\newblock {Perturbations of flexible {L}att\`es maps}.
\newblock {\em Bull. Soc. Math. France} {\bf 141}(2013), 603--614.

\bibitem[CFP]{cfp:fsr}
J.~W. Cannon, W.~J. Floyd, and W.~R. Parry.
\newblock {Finite subdivision rules}.
\newblock {\em Conform. Geom. Dyn.} {\bf 5}(2001), 153--196 (electronic).

\bibitem[CFPP]{cfpp:net}
J.~W. Cannon, W.~J. Floyd, W.~R. Parry, and K.~M. Pilgrim.
\newblock {Nearly {E}uclidean {T}hurston maps}.
\newblock {\em Conform. Geom. Dyn.} {\bf 16}(2012), 209--255.

\bibitem[CGN{\etalchar{+}}]{kmp:fixed}
Kristin Cordwell, Selina Gilbertson, Nicholas Nuechterlein, Kevin~M. Pilgrim,
  and Samantha Pinella.
\newblock {On the classification of critically fixed rational maps}.
\newblock {\em Conform. Geom. Dyn.} {\bf 19}(2015), 51--94.

\bibitem[Cui]{cui-constant-sigma-slides}
Guizhen Cui.
\newblock {Rational maps with (constant) pullback map}.
\newblock
  \url{http://www.math.ac.cn/kyry/cgz/201501/W020150121410977178324.pdf}.

\bibitem[DH]{DH1}
A.~Douady and John Hubbard.
\newblock {A Proof of {T}hurston's Topological Characterization of Rational
  Functions}.
\newblock {\em Acta. Math.} {\bf 171}(1993), 263--297.

\bibitem[FLP]{MR3053012}
Albert Fathi, Fran\c{c}ois Laudenbach, and Valentin Po\'{e}naru.
\newblock {\em Thurston's work on surfaces}, volume~48 of {\em Mathematical
  Notes}.
\newblock Princeton University Press, Princeton, NJ, 2012.
\newblock Translated from the 1979 French original by Djun M. Kim and Dan
  Margalit.

\bibitem[FPP]{Floyd:2018aa}
William Floyd, Walter Parry, and Kevin~M. Pilgrim.
\newblock {Rationality is decidable for nearly Euclidean Thurston maps}.
\newblock {\em to appear, Geom. Dedicata} (2020).

\bibitem[Har]{MR624817}
W.~J. Harvey.
\newblock {Boundary structure of the modular group}.
\newblock In {\em Riemann surfaces and related topics: {P}roceedings of the
  1978 {S}tony {B}rook {C}onference ({S}tate {U}niv. {N}ew {Y}ork, {S}tony
  {B}rook, {N}.{Y}., 1978)}, volume~97 of {\em Ann. of Math. Stud.}, pages
  245--251. Princeton Univ. Press, Princeton, N.J., 1981.

\bibitem[Hlu]{Hlushchanka:2019aa}
Mikhail Hlushchanka.
\newblock {Tischler graphs of critically fixed rational maps and their
  applications, \url{https://arxiv.org/pdf/1904.04759.pdf}}, 04 2019.

\bibitem[KL]{lodge:kelsey:quadratic}
Gregory Kelsey and Russell Lodge.
\newblock {Quadratic Thurston maps with few postcritical points,
  \url{https://arxiv.org/abs/1704.03929}}.
\newblock (04 2017).

\bibitem[KPS]{kmp:kps}
Sarah Koch, Kevin~M. Pilgrim, and Nikita Selinger.
\newblock {Pullback invariants of {T}hurston maps}.
\newblock {\em Trans. Amer. Math. Soc.} {\bf 368}(2016), 4621--4655.

\bibitem[Lod]{lodge:thesis:ecgd}
Russell Lodge.
\newblock {Boundary values of the {T}hurston pullback map}.
\newblock {\em Conform. Geom. Dyn.} {\bf 17}(2013), 77--118.

\bibitem[MBI]{bonk:hlushchanka:iseli:curves}
M.~Hlushchanka M.~Bonk and A.~Iseli.
\newblock {Thurston maps and the dynamics of curves}.
\newblock research announcement, Quasiworld seminar, 2020.

\bibitem[Mal]{saenz:thesis}
Edgar Arturo~Saenz Maldonado.
\newblock {On Nearly Euclidean Thurston Maps}.
\newblock {\em PhD Thesis, Virginia Polytechnic University} (2012).

\bibitem[Mas]{masur:duke:1976}
Howard Masur.
\newblock {Extension of the {W}eil-{P}etersson metric to the boundary of
  {T}eichmuller space}.
\newblock {\em Duke Math. J.} {\bf 43}(1976), 623--635.

\bibitem[Nek]{vn:models:MR3199801}
Volodymyr Nekrashevych.
\newblock {Combinatorial models of expanding dynamical systems}.
\newblock {\em Ergodic Theory Dynam. Systems} {\bf 34}(2014), 938--985.

\bibitem[Par]{parry:netmap:slopefunctions}
Walter Parry.
\newblock {NET map slope functions}.
\newblock September 2018.

\bibitem[Pil1]{kmp:canonical}
Kevin~M. Pilgrim.
\newblock {Canonical {T}hurston obstructions}.
\newblock {\em Adv. Math.} {\bf 158}(2001), 154--168.

\bibitem[Pil2]{kmp:cds}
Kevin~M. Pilgrim.
\newblock {\em Combinations of complex dynamical systems}, volume 1827 of {\em
  Lecture Notes in Mathematics}.
\newblock Springer-Verlag, Berlin, 2003.

\bibitem[Pil3]{kmp:tw:published}
Kevin~M. Pilgrim.
\newblock {An algebraic formulation of {T}hurston's characterization of
  rational functions}.
\newblock {\em Ann. Fac. Sci. Toulouse Math. (6)} {\bf 21}(2012), 1033--1068.

\bibitem[PT]{kmp-tan:blow}
Kevin~M. Pilgrim and Lei Tan.
\newblock {Combining rational maps and controlling obstructions}.
\newblock {\em Ergodic Theory Dynam. Systems} {\bf 18}(1998), 221--245.

\bibitem[Sel]{selinger:pullback}
Nikita Selinger.
\newblock {Thurston's pullback map on the augmented {T}eichm\"uller space and
  applications}.
\newblock {\em Invent. Math.} {\bf 189}(2012), 111--142.

\end{thebibliography}

\end{document}